\documentclass{llncs}
\usepackage{bm,amsmath,amssymb}
\newcommand{\R}{\textbf{R}}

\newcommand{\rank}{\textrm{\upshape rank}}
\newcommand{\diag}{\textrm{diag}}

\spnewtheorem{algorithm}[theorem]{Algorithm}{\bfseries}{\upshape}

\begin{document}

\title{GPGCD, an Iterative Method for Calculating Approximate GCD, for
  Multiple Univariate Polynomials}
\author{Akira Terui}
\institute{Graduate School of Pure and Applied Sciences\\
  University of Tsukuba\\
  Tsukuba, 305-8571, Japan\\
  \email{terui@math.tsukuba.ac.jp}
}
\maketitle

\begin{abstract}
  We present an extension of our GPGCD method, an iterative method for
  calculating approximate greatest common divisor (GCD) of univariate
  polynomials, to multiple polynomial inputs.  For a given pair of
  polynomials and a degree, our algorithm finds a pair of polynomials
  which has a GCD of the given degree and whose coefficients are
  perturbed from those in the original inputs, making the
  perturbations as small as possible, along with the GCD.  In our
  GPGCD method, the problem of approximate GCD is transferred to a
  constrained minimization problem, then solved with the so-called
  modified Newton method, which is a generalization of the
  gradient-projection method, by searching the solution iteratively.
  In this paper, we extend our method to accept more than two
  polynomials with the real coefficients as an input.
\end{abstract}

\section{Introduction}

For algebraic computations on polynomials and matrices, approximate
algebraic algorithms are attracting broad range of attentions
recently.  These algorithms take inputs with some ``noise'' such as
polynomials with floating-point number coefficients with rounding
errors, or more practical errors such as measurement errors, then,
with minimal changes on the inputs, seek a meaningful answer that
reflect desired property of the input, such as a common factor of a
given degree.  By this characteristic, approximate algebraic
algorithms are expected to be applicable to more wide range of
problems, especially those to which exact algebraic algorithms were
not applicable.

As an approximate algebraic algorithm, we consider calculating the
approximate greatest common divisor (GCD) of univariate polynomials,
such that, for a given pair of polynomials and a degree $d$, finding a
pair of polynomials which has a GCD of degree $d$ and whose
coefficients are perturbations from those in the original inputs, with
making the perturbations as small as possible, along with the GCD.
This problem has been extensively studied with various approaches
including the Euclidean method on the polynomial remainder sequence
(PRS) (\cite{bec-lab1998b}, \cite{sas-nod89}, \cite{sch1985}), the
singular value decomposition (SVD) of the Sylvester matrix
(\cite{cor-gia-tra-wat1995}, \cite{emi-gal-lom1997}), the QR
factorization of the Sylvester matrix or its displacements
(\cite{cor-wat-zhi2004}, \cite{zar-ma-fai2000}, \cite{zhi2003}),
Pad\'e approximation \cite{pan2001b}, optimization strategies
(\cite{chi-cor-cor1998}, \cite{kal-yan-zhi2006},
\cite{kal-yan-zhi2007}, \cite{kar-lak1998}, \cite{zen0000a}).
Furthermore, stable methods for ill-conditioned problems have been
discussed (\cite{cor-wat-zhi2004}, \cite{ohs-sug-tor1997},
\cite{san-sas2007}).

Among methods in the above, we focus our attention on optimization
strategies.  Already proposed algorithms utilize iterative methods
including the Levenberg-Marquardt method \cite{chi-cor-cor1998}, the
Gauss-Newton method \cite{zen0000a} and the structured total least
norm (STLN) method (\cite{kal-yan-zhi2006}, \cite{kal-yan-zhi2007}).
Among them, STLN-based methods have shown good performance calculating
approximate GCD with sufficiently small perturbations efficiently.

In this paper, we discuss an extension of the GPGCD method, proposed
by the present author (\cite{ter2009}, \cite{ter2009b}), an iterative
method with transferring the original approximate GCD problem into a
constrained optimization problem, then solving it by the so-called
modified Newton method \cite{tan1980}, which is a generalization of
the gradient-projection method \cite{ros1961}.  In the previous papers
(\cite{ter2009}, \cite{ter2009b}), we have shown that our method
calculates approximate GCD with perturbations as small as those
calculated by the STLN-based methods and with significantly better
efficiency than theirs.  While our previous methods accept two
polynomials with the real or the complex coefficients as inputs and
outputs, respectively, we extend it to handle more than two polynomial
inputs with the real coefficients in this paper.

The rest part of the paper is organized as follows.  In
Section~\ref{sec:formulation}, we transform the approximate GCD
problem into a constrained minimization problem for the case with the
complex coefficients.  In Section~\ref{sec:gpgcd}, we show details for
calculating the approximate GCD, with discussing issues in
minimizations.  In Section~\ref{sec:exp}, we demonstrate performance
of our algorithm with experiments.

\section{Formulation of the Approximate GCD Problem}
\label{sec:formulation}

Let $P_1(x),\ldots,P_n(x)$ be real univariate polynomials of degree
$d_1,\ldots,d_n$, respectively, given as
\begin{equation*}
    P_i(x) = p_{d_i}^{(i)}x^{d_i} + \cdots p_1^{(i)} x + p_0^{(i)},
\end{equation*}
for $i=1,\ldots,n$, with $\min\{d_1,\ldots,d_n\}>0$.  We permit $P_i$
and $P_j$ be relatively prime for any $i\ne j$ in general.  For a
given integer $d$ satisfying $\min\{d_1,\ldots,d_n\}>d>0$, let us
calculate a deformation of $P_1(x),\ldots,P_n(x)$ in the form of
\begin{equation*}
    \tilde{P_i}(x) = P_i(x) + \varDelta P_i(x) = H(x)\cdot
    \bar{P_i}(x), 
\end{equation*}
where $\varDelta P_i(x)$ is a real polynomial whose degrees do not
exceed $d_i$, respectively, $H(x)$ is a polynomial of degree $d$, and
$\bar{P_i}(x)$ and $\bar{P_j}(x)$ are pairwise relatively prime for
any $i\ne j$.  In this situation, $H(x)$ is an approximate GCD of
$P_1(x),\ldots,P_n(x)$.  For a given $d$, we try to minimize
$\|\varDelta P_1(x)\|_2^2+\cdots+\|\varDelta P_n(x)\|_2^2$, the norm
of the deformations.

For a real univariate polynomial $P(x)$ represented as
$
P(x) = p_nx^n+\cdots+p_0x^0,
$
let $C_k(P)$ be a real $(n+k,k+1)$ matrix defined as
\[
\begin{array}{ccl}
  C_k(P) & = & 
  \begin{pmatrix}
    p_n & & \\
    \vdots & \ddots & \\
    p_0 & & p_n \\
    & \ddots & \vdots \\
    & & p_0
  \end{pmatrix}
  ,
  \\[-2mm]
  & & 
  \hspace{3mm}
  \underbrace{\hspace{13.5mm}}_{k+1}
\end{array}
\]
and let $\bm{p}$ be the coefficient vector of $P(x)$ defined as
\begin{equation}
  \label{eq:coefvecdef}
  \bm{p}=(p_n,\ldots,p_0).
\end{equation}

In this paper, for a generalized Sylvester matrix,  we use a
formulation by Rupprecht \cite[Sect.~3]{rup1999}.  Then, a generalized
Sylvester matrix for $P_1,\ldots,P_n$ becomes as
\begin{equation}
  \label{eq:sylvestermat-n}
  N(P_1,\ldots,P_n)=
  \begin{pmatrix}
    C_{d_1-1}(P_2) & C_{d_2-1}(P_1) & \bm{0} & \cdots & \bm{0} \\
    C_{d_1-1}(P_3) & \bm{0} & C_{d_3-1}(P_1) & \cdots & \bm{0} \\
    \vdots & \vdots & & \ddots & \vdots \\
    C_{d_1-1}(P_n) & \bm{0} & \cdots & \bm{0} & C_{d_n-1}(P_1)
  \end{pmatrix}
  ,
\end{equation}
and the $k$-th subresultant matrix (with $\min\{d_1,\ldots,d_n\}>k\ge
0$) is also defined similarly as 
\begin{multline}
  \label{eq:subresmat-n}
  N_k(P_1,\ldots,P_n)
  \\
  =
  \begin{pmatrix}
    C_{d_1-1-k}(P_2) & C_{d_2-1-k}(P_1) & \bm{0} & \cdots & \bm{0} \\
    C_{d_1-1-k}(P_3) & \bm{0} & C_{d_3-1-k}(P_1) & \cdots & \bm{0} \\
    \vdots & \vdots & & \ddots & \vdots \\
    C_{d_1-1-k}(P_n) & \bm{0} & \cdots & \bm{0} & C_{d_n-1-k}(P_1)
  \end{pmatrix}
  ,
\end{multline}
with
\begin{equation}
  \label{eq:rk}
  r_k=d_1+d_2+\cdots+d_n-(n-1)k+(n-2)d_1
\end{equation}
rows and
\begin{equation}
  \label{eq:ck}
  c_k=d_1+d_2+\cdots+d_n-n\cdot k 
\end{equation}
columns. 

Calculation of GCD is based on the following fact.
\begin{proposition}[Rupprecht {\cite[Proposition 3.1]{rup1999}}]
  \label{prop:deggcd}
  $N_k(P_1,\ldots,P_n)$ has full rank if and only if
  $\deg(\gcd(P_1,\ldots,P_n))\le k$.
\end{proposition}

Thus, for a given degree $d$, if
$N_{d-1}(\tilde{P_1},\ldots,\tilde{P_n})$ is rank-deficient, then
there exist real univariate polynomials $U_1(x),\ldots,U_n(x)$ of
degree at most $d_1-d,\ldots,d_n-d$, respectively, satisfying 
\begin{equation}
  \label{eq:uicond}
  U_1\tilde{P}_i+U_i\tilde{P}_1=0,
\end{equation}
for $i=2,\ldots,n$.  In such a case, if $U_i$ and $U_j$ are pairwise
relatively prime for any $i\ne j$, then
$H=\frac{\tilde{P}_1}{U_1}=-\frac{\tilde{P}_2}{U_2}=\cdots=
-\frac{\tilde{P}_n}{U_n}$ becomes the expected GCD.  Therefore, for
given polynomials $P_1,\ldots,P_n$ and a degree $d$, our problem is to
find perturbations $\varDelta P_1,\ldots,\varDelta P_n$ along with
cofactors $U_1,\ldots,U_n$ satisfying \eqref{eq:uicond} with making
$\|\varDelta P_1(x)\|_2^2+\cdots+\|\varDelta P_n(x)\|_2^2$ as small as
possible.

By representing $\tilde{P}_i(x)$ and $U_i(x)$ as
\begin{equation}
  \label{eq:tildepiui}
  \begin{split}
    \tilde{P}_i(x) &
    =\tilde{p}_{d_i}^{(i)}x^{d_i}+\cdots+\tilde{p}_1^{(i)}x+\tilde{p}_0^{(i)}, 
    \\
    U_i(x) & =u_{d_i-d}^{(i)}x^{d_i-d}+\cdots+u_1^{(i)}x+u_0^{(i)},
  \end{split}
\end{equation}
we express the objective function and the constraint as follows.  For
the objective function, $\|\varDelta P_1(x)\|_2^2+\cdots+\|\varDelta
P_n(x)\|_2^2$ becomes as
\begin{equation}
  \label{eq:objective-multiple-real}
  \|\varDelta P_1(x)\|_2^2+\cdots+\|\varDelta P_n(x)\|_2^2
  =
  \sum_{i=1}^n
  \left\{
    \sum_{j=0}^{d_i}
    \left(
      \tilde{p}_j^{(i)}-p_j^{(i)}
    \right)^2
  \right\}
  .
\end{equation}
For the constraint, \eqref{eq:uicond} becomes as
\begin{equation}
  \label{eq:ucond-real}
  N_{d-1}(\tilde{P}_1,\ldots,\tilde{P}_n)\cdot{}^t
  (\bm{u}_1,\ldots,\bm{u}_n)=\bm{0},
\end{equation}
where $\bm{u}_i$ is the coefficient vector of $U_i(x)$ defined as in
\eqref{eq:coefvecdef}.  Furthermore, we add another constraint for the
coefficient of $U_i(x)$ as
\begin{equation}
  \label{eq:uiconstraint-real}
  \|U_1\|_2^2+\cdots+\|U_n\|_2^2=1,
\end{equation}
which can be represented together with \eqref{eq:ucond-real} as
\begin{equation}
  \label{eq:ucond-2-real}
  \begin{pmatrix}
    \bm{u}_1\quad \cdots\quad \bm{u_n} & -1 \\
    N_{d-1}(\tilde{P}_1,\ldots,\tilde{P}_n) & \bm{0}
  \end{pmatrix}
  \cdot{}^t
  (\bm{u}_1,\ldots,\bm{u}_n,1)=\bm{0},
\end{equation}
where \eqref{eq:uiconstraint-real} has been put on the top of
\eqref{eq:ucond-real}.  Note that, in \eqref{eq:ucond-2-real}, we have
total of
\begin{equation}
  \label{eq:numconstraint}
  \bar{d}=d_1+\cdots+d_n-(n-1)(d-1)+(n-2)d_1+1
\end{equation}
equations in the coefficients of polynomials in
\eqref{eq:tildepiui} as a constraint, with the $j$-th row of which
is expressed as $g_j=0$.

Now, we substitute the variables
\begin{equation}
  \label{eq:var0-real}
  (
  \tilde{p}_{d_1}^{(1)},\ldots,\tilde{p}_0^{(1)},
  \ldots,
  \tilde{p}_{d_n}^{(n)},\ldots,\tilde{p}_0^{(n)},
  u_{d_1-d}^{(1)},\ldots,u_0^{(1)},
  \ldots,
  u_{d_n-d}^{(n)},\ldots,u_0^{(n)}
  ),
\end{equation}
as $\bm{x}=(x_1,\ldots,x_{2(d_1+\cdots+d_n)+(2-d)n})$, then
\eqref{eq:objective-multiple-real} and \eqref{eq:ucond-2-real} become
as
\begin{multline}
  \label{eq:objective-multiple-real-2}
  f(\bm{x})=
  (x_1-p_{d_1}^{(1)})^2+\cdots+(x_{d_1}-p_0^{(1)})^2
  +\cdots
  \\
  \cdots+
  (x_{d_1+\cdots+d_{n-1}+n}-p_{d_n}^{(n)})^2+\cdots
  +(x_{d_1+\cdots+d_{n-1}+d_n+n}-p_0^{(n)})^2,
\end{multline}
\begin{equation}
  \label{eq:constraint-multiple-real}
  \bm{g}(\bm{x})=
  {}^t(g_1(\bm{x}), \ldots, g_{\bar{d}}(\bm{x}))
  =
  \bm{0},
\end{equation}
respectively, where $\bar{d}$ in \eqref{eq:constraint-multiple-real} is
defined as in \eqref{eq:numconstraint}.  Therefore, the problem of
finding an approximate GCD can be formulated as a constrained
minimization problem of finding a minimizer of the objective function
$f(\bm{x})$ in (\ref{eq:objective-multiple-real-2}), subject to
$\bm{g}(\bm{x})=\bm{0}$ in (\ref{eq:constraint-multiple-real}).

\section{The Algorithm for Approximate GCD}
\label{sec:gpgcd}

We calculate an approximate GCD by solving the constrained
minimization problem \eqref{eq:objective-multiple-real-2},
\eqref{eq:constraint-multiple-real} with the gradient projection
method by Rosen \cite{ros1961} (whose initials become the name of our
GPGCD method) or the modified Newton method by Tanabe \cite{tan1980}
(for review, see the author's previous paper \cite{ter2009}).  Our
preceding experiments (\cite[Sect.\ 5.1]{ter2009},
\cite[Sect.\ 4]{ter2009b}) have shown that the modified Newton method was more
efficient than the original gradient projection method while the both
methods have shown almost the same convergence property, thus we adopt
the modified Newton method in this paper.

In applying the modified Newton method to the approximate GCD problem,
we discuss issues in the construction of the algorithm in detail, such
as
\begin{itemize}
\item Representation of the Jacobian matrix $J_{\bm{g}}(\bm{x})$ and
  certifying that $J_{\bm{g}}(\bm{x})$ has full rank
  (Sect.~\ref{sec:jacobianmat}),
\item Setting the initial values (Sect.~\ref{sec:init}),
\item Regarding the minimization problem as the minimum distance
  problem (Sect.~\ref{sec:mindistance}),
\item Calculating the actual GCD and correcting the coefficients of
  $\tilde{P}_i$ (Sect.~\ref{sec:correct}),
\end{itemize}
as follows.

\subsection{Representation and the rank of the  Jacobian Matrix}
\label{sec:jacobianmat}

By the definition of the constraint
\eqref{eq:constraint-multiple-real}, we have the Jacobian matrix
$J_{\bm{g}}(\bm{x})$ (with the original notation of variables
\eqref{eq:var0-real} for $\bm{x}$) as
\begin{multline*}
  J_{\bm{g}}(\bm{x}) =
  \left(
    \begin{array}{ccccc}
      \bm{0} & \bm{0} & \bm{0} & \cdots & \bm{0} 
      \\ 
      C_{d_1}(U_2) & C_{d_2}(U_1) & \bm{0} & \cdots &  \bm{0} \\
      C_{d_1}(U_3) & \bm{0} & C_{d_3}(U_1) & &  \bm{0} \\
      \vdots & \vdots  & & \ddots & \vdots \\
      C_{d_1}(U_n) & \bm{0} & \cdots & \bm{0} & C_{d_n}(U_1)
    \end{array}
  \right.
  \\
  \left.
    \begin{array}{ccccc}
      2\cdot{}^t\bm{u}_1 & 2\cdot{}^t\bm{u}_2 & 2\cdot{}^t\bm{u}_3 &
      \cdots & 2\cdot{}^t\bm{u}_n
      \\ 
      C_{d_1-d}(P_2) & C_{d_2-d}(P_1) & \bm{0} & \cdots &  \bm{0} \\
      C_{d_1-d}(P_3) & \bm{0} & C_{d_3-d}(P_1) & &  \bm{0} \\
      \vdots & \vdots  & & \ddots & \vdots \\
      C_{d_1-d}(P_n) & \bm{0} & \cdots & \bm{0} & C_{d_n-d}(P_1)
    \end{array}
  \right)
  ,
\end{multline*}
which can easily be constructed in every iteration. Note that the
number of rows in $J_{\bm{g}}(\bm{x})$ is equal to $\bar{d}$ in
\eqref{eq:numconstraint}, which is equal to the number of constraints,
while the number of columns is equal to $2(d_1+\cdots+d_n)+(2-d)n$,
which is equal to the number of variables (see \eqref{eq:var0-real}).

In executing iterations, we need to keep that $J_{\bm{g}}(\bm{x})$ has
full rank: otherwise, we are unable to decide proper search direction.
For this requirement, we have the following observations.
\begin{proposition}
  \label{prop:fullrank}
  Assume that we have $\deg d<\min\{d_1,\ldots,d_n\}-1$ and $\deg
  U_i\ge 1$ for $i=1,\ldots,n$.  Let $\bm{x}^*\in V_{\bm{g}}$ be any
  feasible point satisfying \eqref{eq:constraint-multiple-real}.
  Then, if the corresponding polynomials do not have a GCD whose
  degree exceeds $d$, then $J_{\bm{g}}(\bm{x}^*)$ has full rank.
\end{proposition}
\begin{proof}
  Let
  \[
  \bm{x}^*=
  (
  \tilde{p}_{d_1}^{(1)},\ldots,\tilde{p}_0^{(1)},
  \ldots,
  \tilde{p}_{d_n}^{(n)},\ldots,\tilde{p}_0^{(n)},
  u_{d_1-d}^{(1)},\ldots,u_0^{(1)},
  \ldots,
  u_{d_n-d}^{(n)},\ldots,u_0^{(n)}
  )
  \]
  as in \eqref{eq:var0-real}, with its polynomial representation
  expressed as in \eqref{eq:tildepiui} (note that this assumption
  permits the polynomials $\tilde{P}_i(x)$ to be relatively prime in
  general).  To verify our claim, we show that we have
  $\rank(J_{\bm{g}}(\bm{x}^*))=\bar{d}=d_1+\cdots+d_n-(n-1)(d-1)+(n-2)d_1+1$
  (see \eqref{eq:numconstraint}).  Let us divide
  $J_{\bm{g}}(\bm{x}^*)$ into two column blocks such that
  $J_{\bm{g}}(\bm{x}^*)=
  \begin{pmatrix}
    J_\mathrm{L} \mid J_\mathrm{R}
  \end{pmatrix}
  $, where $J_\mathrm{L}$ and $J_\mathrm{R}$ are expressed as
  \[
  \begin{split}
    J_\mathrm{L} &=
    \begin{pmatrix}
      \bm{0} & \bm{0} & \bm{0} & \cdots & \bm{0} 
      \\ 
      C_{d_1}(U_2) & C_{d_2}(U_1) & \bm{0} & \cdots &  \bm{0} \\
      C_{d_1}(U_3) & \bm{0} & C_{d_3}(U_1) & &  \bm{0} \\
      \vdots & \vdots  & & \ddots & \vdots \\
      C_{d_1}(U_n) & \bm{0} & \cdots & \bm{0} & C_{d_n}(U_1)
    \end{pmatrix}
    ,
    \\
    J_\mathrm{R} &=
    \begin{pmatrix}
      2\cdot{}^t\bm{u}_1 & 2\cdot{}^t\bm{u}_2 & 2\cdot{}^t\bm{u}_3 &
      \cdots & 2\cdot{}^t\bm{u}_n
      \\ 
      C_{d_1-d}(P_2) & C_{d_2-d}(P_1) & \bm{0} & \cdots &  \bm{0} \\
      C_{d_1-d}(P_3) & \bm{0} & C_{d_3-d}(P_1) & &  \bm{0} \\
      \vdots & \vdots  & & \ddots & \vdots \\
      C_{d_1-d}(P_n) & \bm{0} & \cdots & \bm{0} & C_{d_n-d}(P_1)
    \end{pmatrix}
    ,
  \end{split}
  \]
  respectively.  Then, we have the following lemma.
  \begin{lemma}
    \label{lem:fullrank}
    We have
    $\rank(J_\mathrm{L})=\bar{d}=d_1+\cdots+d_n-(n-1)(d-1)+(n-2)d_1$. 
  \end{lemma}
  \begin{proof}
    Let $\bar{J}$ be a submatrix of $J_\mathrm{L}$ by eliminating the
    top row.  Since the number of rows in $J_\mathrm{L}$ is equal to 
    $\bar{d}=d_1+\cdots+d_n-(n-1)(d-1)+(n-2)d_1$, we show that
    $\bar{J}$ has full rank.

    For $i=2,\ldots,n$, let us divide column blocks $C_{d_1}(U_i)$ and
    $C_{d_i}(U_1)$ as
    \begin{align}
      C_{d_1}(U_i) &=
       \begin{array}[b]{l}
         \hspace{2mm}
         \substack{
           d+1\\ 
           \overbrace{\hspace{13mm}}
         }
         \hspace{2mm}
         \substack{
           d_1-d
           \\
           \overbrace{\hspace{12.5mm}}
         }           
         \\
         \begin{pmatrix}
           C_{d_1}(U_i)_\mathrm{L} &  C_{d_1}(U_i)_\mathrm{R}           
         \end{pmatrix}
       \end{array}
       ,
       \nonumber\\
       C_{d_1}(U_i)_\mathrm{L} &=
       \begin{pmatrix}
         C_{d}(U_i)
         \\[1mm]
         \bm{0}
       \end{pmatrix}
       \begin{matrix}
         \\[1mm]
         \} \substack{d_1-d}
       \end{matrix}
       \;
       ,
       \quad
       C_{d_1}(U_i)_\mathrm{R} =
       \begin{pmatrix}
         \bm{0}
         \\[1mm]
         C_{d_1-d-1}(U_i)
       \end{pmatrix}
       \begin{matrix}
         \} \substack{d+1}
         \\[1mm]
         \} \substack{d_1+d_i-2d}
       \end{matrix}
       \;
       ,
       \label{eq:cd1uilr}\\
      C_{d_i}(U_1) &=
       \begin{array}[b]{l}
         \hspace{2mm}
         \substack{
           d+1\\ 
           \overbrace{\hspace{13mm}}
         }
         \hspace{2mm}
         \substack{
           d_1-d
           \\
           \overbrace{\hspace{12.5mm}}
         }           
         \\
         \begin{pmatrix}
           C_{d_i}(U_1)_\mathrm{L} &  C_{d_i}(U_1)_\mathrm{R}           
         \end{pmatrix}
       \end{array}
       ,
       \nonumber\\
       C_{d_i}(U_1)_\mathrm{L} &=
       \begin{pmatrix}
         C_{d}(U_1)
         \\[1mm]
         \bm{0}
       \end{pmatrix}
       \begin{matrix}
         \\[1mm]
         \} \substack{d_1-d}
       \end{matrix}
       \;
       ,
       \quad
       C_{d_i}(U_1)_\mathrm{R} =
       \begin{pmatrix}
         \bm{0}
         \\[1mm]
         C_{d_i-d-1}(U_1)
       \end{pmatrix}
       \begin{matrix}
         \} \substack{d+1}
         \\[1mm]
         \} \substack{d_1+d_i-2d}
       \end{matrix}
       \;
       \label{eq:cdiu1lr}
       ,
    \end{align}
    respectively, thus $\bar{J}$ is expressed as
    \begin{multline*}
      \bar{J}=
      \left(
        \begin{array}{cc|cc|ccc}
          C_{d_1}(U_2)_\mathrm{L} & C_{d_1}(U_2)_\mathrm{R} &
          C_{d_2}(U_1)_\mathrm{L} & C_{d_2}(U_1)_\mathrm{R} &
          \bm{0} & \bm{0} & \cdots \\  
          C_{d_1}(U_3)_\mathrm{L} & C_{d_1}(U_3)_\mathrm{R} &
          \bm{0} & \bm{0} & C_{d_3}(U_1)_\mathrm{L} &
          C_{d_3}(U_1)_\mathrm{L} & \\
          \vdots & \vdots  & \vdots & \vdots  & & \ddots & \ddots \\ 
          C_{d_1}(U_n)_\mathrm{L} & C_{d_1}(U_n)_\mathrm{R} &
          \bm{0} & \bm{0} & \bm{0} & \bm{0} & \cdots
        \end{array}
      \right.
      \\
      \left.
        \begin{array}{ccc|cc}
          \cdots & \bm{0} & \bm{0} & \bm{0} &  \bm{0} \\
          \ddots & \ddots &  & \vdots & \vdots \\
          & C_{d_{n-1}}(U_1)_\mathrm{L} &
          C_{d_{n-1}}(U_1)_\mathrm{R} & \bm{0} & \bm{0} \\ 
          \cdots & \bm{0} & \bm{0} & C_{d_n}(U_1)_\mathrm{R} &
          C_{d_n}(U_1)_\mathrm{R} 
        \end{array}
      \right)
      .
    \end{multline*}
    Then, by exchanges of columns, we can transform $\bar{J}$ to
    $
    \hat{J}=
    \begin{pmatrix}
      \hat{J}_\mathrm{L} & \hat{J}_\mathrm{R}
    \end{pmatrix}
    $,
    where
    \begin{align*}
      \hat{J}_\mathrm{L} &=
      \begin{pmatrix}
        C_{d_1}(U_2)_\mathrm{L} & C_{d_2}(U_1)_\mathrm{L} & \bm{0} &
        \cdots & \bm{0} \\ 
        C_{d_1}(U_3)_\mathrm{L} & \bm{0} & C_{d_3}(U_1)_\mathrm{L} & &
        \bm{0} \\ 
        \vdots & \vdots  & & \ddots & \vdots \\
        C_{d_1}(U_n)_\mathrm{L} & \bm{0} & \cdots & \bm{0} &
        C_{d_n}(U_1)_\mathrm{L} 
      \end{pmatrix}
      ,
      \\
      \hat{J}_\mathrm{R} &=
      \begin{pmatrix}
        C_{d_1}(U_2)_\mathrm{R} & C_{d_2}(U_1)_\mathrm{R} & \bm{0} &
        \cdots & \bm{0} \\ 
        C_{d_1}(U_3)_\mathrm{R} & \bm{0} & C_{d_3}(U_1)_\mathrm{R} & &
        \bm{0} \\ 
        \vdots & \vdots  & & \ddots & \vdots \\
        C_{d_1}(U_n)_\mathrm{R} & \bm{0} & \cdots & \bm{0} &
        C_{d_n}(U_1)_\mathrm{R} 
      \end{pmatrix}
      .
    \end{align*}

    We see that nonempty rows in $\hat{J}_\mathrm{R}$ consist of
    $N(U_1,\ldots,U_n)$, a generalized Sylvester matrix for
    $U_1,\ldots,U_n$ (see \eqref{eq:sylvestermat-n}).  By the
    assumption, $U_1,\ldots,U_n$ are pairwise relatively prime, thus,
    by Prop.~\ref{prop:deggcd}, $\rank(\hat{J}_\mathrm{R})$ is equal
    to the number of nonempty rows in $\hat{J}_\mathrm{R}$, which is
    equal to $d_2+\cdots+d_n+(n-1)(d_1-2d)$ (see \eqref{eq:cd1uilr}
    and \eqref{eq:cdiu1lr}).

    On the other hand, in $\hat{J}_L$, column blocks
    $C_{d_2}(U_1)_\mathrm{L},C_{d_3}(U_1)_\mathrm{L},\ldots,C_{d_n}(U_1)_\mathrm{L}$
    are lower triangular matrices with $d+1$ diagonal elements, which
    shows that $\rank(\hat{J}_L)$ is equal to the sum of
    the number of columns in
    $C_{d_2}(U_1)_\mathrm{L},C_{d_3}(U_1)_\mathrm{L},\ldots$,
    $C_{d_n}(U_1)_\mathrm{L}$, which is equal to $(n-1)(d+1)$.

    Furthermore, we see that the row position of diagonal elements in
    $C_{d_2}(U_1)_\mathrm{L}$,
    $C_{d_3}(U_1)_\mathrm{L},\ldots,C_{d_n}(U_1)_\mathrm{L}$
    correspond to the position of the empty rows in
    $\hat{J}_\mathrm{R}$, thus the columns in
    $C_{d_2}(U_1)_\mathrm{L},C_{d_3}(U_1)_\mathrm{L},\ldots,
    C_{d_n}(U_1)_\mathrm{L}$ are linearly independent along with the
    columns in $\hat{J}_R$.  Therefore, we have
    \[
    \rank(\bar{J})=\rank(\hat{J}_L)+\rank(\hat{J}_R)
    =
    d_1+\cdots+d_n-(n-1)(d-1)+(n-2)d_1,
    \]
    which proves the lemma.
  \end{proof}

  \textsc{Proof of Proposition~\ref{prop:fullrank} (continued).}  By
  the assumptions, we have at least one nonzero coordinate in the top
  row in $J_\mathrm{R}$, while we have no nonzero coordinate in the
  top row in $J_\mathrm{L}$, thus we have
  $\rank(J_{\bm{g}}(\bm{x}))=d_1+\cdots+d_n-(n-1)(d-1)+(n-2)d_1+1$,
  which proves the proposition.  \qed
\end{proof}

Proposition~\ref{prop:fullrank} says that, under certain conditions,
so long as the search direction in the minimization problem satisfies
that corresponding polynomials have a GCD of degree not exceeding $d$,
then $J_{\bm{g}}(\bm{x})$ has full rank, thus we can safely calculate
the next search direction for approximate GCD.

\subsection{Setting the Initial Values}
\label{sec:init}

At the beginning of iterations, we give the initial value $\bm{x}_0$
by using the singular value decomposition (SVD) \cite{dem1997} of
$N_{d-1}(P_1,\ldots,P_n)$ (see \eqref{eq:subresmat-n}) as
$
    N_{d-1} =  U\,\Sigma\,{}^tV,
    U = (\bm{w}_1,\ldots,\bm{w}_{c_{d-1}}),
$
$
    \Sigma = \diag(\sigma_1,\ldots,\sigma_{c_{d-1}}),
    V = (\bm{v}_1,\ldots,\bm{v}_{c_{d-1}}),
    $ where $\bm{w}_j\in\R^{r_{d-1}}$, $\bm{v}_j\in\R^{c_{d-1}}$ with
    $r_k$ and $c_k$ as in \eqref{eq:rk} and \eqref{eq:ck},
    respectively, and $\Sigma=\diag(\sigma_1,\ldots,$
    $\sigma_{c_{d-1}})$ denotes the diagonal matrix with the $j$-th
    diagonal element of which is $\sigma_j$.  Note that $U$ and $V$
    are orthogonal matrices.  Then, by a property of the SVD
    \cite[Theorem~3.3]{dem1997}, the smallest singular value
    $\sigma_{c_{d-1}}$ gives the minimum distance of the image of the
    unit sphere $\textrm{S}^{c_{d-1}-1}$, given as
$
\textrm{S}^{c_{d-1}-1}=
\{
\bm{x}\in\R^{c_{d-1}} \mid \|\bm{x}\|_2=1
\}
,
$
by $N_{d-1}(P_1,\ldots,P_n)$, represented as
$
N_{d-1}\cdot\textrm{S}^{c_{d-1}-1}=
\{
N_{d-1}\bm{x}\mid \bm{x}\in\R^{c_{d-1}}, \|\bm{x}\|_2=1
\},
$
from the origin, along with $\sigma_{c_{d-1}}\bm{w}_{c_{d-1}}$ as its
coordinates.  Thus, we have
$
N_{d-1}\cdot \bm{v}_{c_{d-1}} = \sigma_{c_{d-1}}\bm{w}_{c_{d-1}}.
$
For
$\bm{v}_{c_{d-1}}=
{}^t(\bar{u}_{d_1-d}^{(1)},\ldots,\bar{u}_0^{(1)},
\ldots,
\bar{u}_{d_n-d}^{(n)},\ldots,\bar{u}_0^{(n)})$, 
let
$\bar{U}_i(x) =
\bar{u}_{d_i-d}^{(i)}x^{d_i-d}+\cdots+\bar{u}_0^{(i)}x^0$
for $i=1,\ldots,n$.
Then, $\bar{U}_1(x),\ldots,\bar{U}_n(x)$ give the least norm of
$U_1P_i+U_iP_1$ satisfying $\|U_1\|_2^2+\cdots+\|U_n\|_2^2=1$ by
putting $U_i(x)=\bar{U}_i(x)$ in \eqref{eq:tildepiui}.

Therefore, we admit the coefficients of $P_1,\ldots,P_n$, $\bar{U}_1,\ldots,\bar{U}_n$ as the initial values of the iterations as
\begin{multline*}
  \bm{x}_0 = 
  (
  p_{d_1}^{(1)},\ldots,p_0^{(1)},
  \ldots,
  p_{d_n}^{(n)},\ldots,p_0^{(n)},
  \bar{u}_{d_1-d}^{(1)},\ldots,\bar{u}_0^{(1)},
  \ldots,
  \bar{u}_{d_n-d}^{(n)},\ldots,\bar{u}_0^{(n)}
  ).
\end{multline*}

\subsection{Regarding the Minimization Problem as the Minimum Distance
  (Least Squares) Problem} 
\label{sec:mindistance}

Since we have the object function $f$ as in
(\ref{eq:objective-multiple-real-2}), we have
\begin{multline*}
  \nabla f(\bm{x})=
  2\cdot {}^t(
  x_1-p_{d_1}^{(1)},\ldots,x_{d_1}-p_0^{(1)},
  \ldots,
  \\
  x_{d_1+\cdots+d_{n-1}+n}-p_{d_n}^{(n)},\ldots,
  x_{d_1+\cdots+d_{n-1}+d_n+n}-p_0^{(n)},
  0,\ldots,0).
\end{multline*}
However, we can regard our problem as finding a point $\bm{x}\in
V_{\bm{g}}$ which has the minimum distance to the initial point
$\bm{x}_0$ with respect to the
$(x_1,\ldots,x_{d_1+\cdots+d_{n-1}+d_n+n})$-coordinates which
correspond to the coefficients in $P_i(x)$.
Therefore, as in the case for two polynomials (see the author's
previous papers (\cite{ter2009}, \cite{ter2009b})), we change the
objective function as $\bar{f}(\bm{x})=\frac{1}{2}f(\bm{x})$, then
solve the minimization problem of $\bar{f}(\bm{x})$, subject to
$\bm{g}(\bm{x})=\bm{0}$.

\subsection{Calculating the Actual GCD and Correcting the Deformed
  Polynomials}
\label{sec:correct}

After successful end of the iterations, we obtain the coefficients of
$\tilde{P}_i(x)$ and $U_i(x)$ satisfying (\ref{eq:uicond}) with
$U_i(x)$ are relatively prime.  Then, we need to compute the actual
GCD $H(x)$ of $\tilde{P}_i(x)$.  Although $H$ can be calculated as the
quotient of $\tilde{P}_i$ divided by $U_i$, naive polynomial division
may cause numerical errors in the coefficient.  Thus, we calculate the
coefficients of $H$ by the so-called least squares division
\cite{zen0000a}, followed by correcting the coefficients in
$\tilde{P}_i$ by using the calculated $H$, as follows.

For polynomials $\tilde{P}_i$, and $U_i$
represented as in \eqref{eq:tildepiui} and $H$ represented as
\[
H(x)=h_d x^d+\cdots+h_0 x^0,
\]
solve the equations $HU_i=\tilde{P}_i$ with respect to $H$ as solving
the least squares problems of a linear system
\begin{equation}
  \label{eq:hsystemi}
  C_d(U_i)\,{}^t(h_d\ldots,h_0) =
  {}^t(p_{d_i}^{(i)},\cdots,p_0^{(i)}).
\end{equation}
Let $H_i(x)\in\R[x]$ be a candidate for the GCD whose coefficients are
calculated as the least squares solutions of \eqref{eq:hsystemi}.
Then, for $i=2,\ldots,n$, calculate the norms of the residues as
\[
r_i = \sum_{j=1}^n\|P_j-H_iU_j\|_2^2,
\]
and set the GCD $H(x)$ be $H_i(x)$ giving the minimum value of $r_i$
so that the perturbed polynomials make the minimum amount of
perturbations in total.

Finally, for the chosen $H(x)$, correct the coefficients of
$\tilde{P}_i(x)$ as $\tilde{P}_i(x)=H(x)\cdot U_i(x)$ for
$i=1,\ldots,n$.

\section{Experiments}
\label{sec:exp}

We have implemented our GPGCD method on the computer algebra system
Maple and compared its performance with a method based on the
structured total least norm (STLN) method \cite{kal-yan-zhi2006} for
randomly generated polynomials with approximate GCD.  The tests have
been carried out on Intel Core2 Duo Mobile Processor T7400 (in Apple
MacBook ``Mid-2007'' model) at $2.16$ GHz with RAM 2GB, under MacOS X
10.5.

In the tests, we have generated random polynomials with GCD then added
noise, as follows.  First, we have generated a monic polynomial
$P_0(x)$ of degree $m$ with the GCD of degree $d$.  The GCD and the
prime parts of degree $m-d$ are generated as monic polynomials and
with random coefficients $c\in[-10,10]$ of floating-point numbers.
For noise, we have generated a polynomial $P_{\mathrm{N}}(x)$ of
degree $m-1$ with random coefficients as the same as for $P_0(x)$.
Then, we have defined a test polynomial $P(x)$ as
$
  P(x) 
  = P_0(x)+\frac{e_P}{\|P_{\mathrm{N}}(x)\|_2}P_{\mathrm{N}}(x),
$
scaling the noise such that the $2$-norm of the noise for $P$ is equal
to $e_P$.  In the present test, we set $e_P=0.1$.  

In this test, we have compared our implementation against a method
based on the structured total least norm (STLN) method
\cite{kal-yan-zhi2006}, using their implementation (see
Acknowledgments).  In their STLN-based method, we have used the
procedure \verb|R_con_mulpoly| which calculates the approximate GCD of
several polynomials in $\R[x]$.  The tests have been carried out on
Maple 13 with \verb|Digits=15| executing hardware floating-point
arithmetic.  For every example, we have generated $10$ random test
cases as in the above.  In executing the GPGCD method, we set $u=100$
and a threshold of the $2$-norm of the ``update'' vector in each
iteration $\varepsilon=1.0\times 10^{-8}$; in \verb|R_con_mulpoly|, we
set the tolerance $e=1.0\times 10^{-8}$.

Table~\ref{tab:appgcd} shows the results of the test.  In each test,
we have given several polynomials of the same degree as the input.
The second column with $(m,d,n)$ denotes the degree of input
polynomials, degree of GCD, and the number of input polynomials,
respectively.  The columns with ``STLN'' are the data for the
STLN-based method, while those with ``GPGCD'' are the data for the
GPGCD method.  ``\#Fail'' is the number of ``failed'' cases such as:
in the STLN-based method, the number of iterations exceeds $50$ times
(which is the built-in threshold in the program), while, in the GPGCD
method, the number of iterations exceeds $100$ times.  All the other
data are the average over results for the ``not failed'' cases:
``Error'' is the sum of perturbation
$\sum_{i=1}^n\|\tilde{P}_i-P_i\|_2^2$, where ``$ae\!-\!b$'' denotes
$a\times 10^{-b}$; ``\#Iterations'' is the number of iterations;
``Time'' is computing time in seconds.

We see that, in the most of tests, both methods calculate approximate
GCD with almost the same amount of perturbations. In the most of
tests, the GPGCD method runs faster than STLN-based method.  However,
running time of the GPGCD method increases as much as that of the
STLN-based method in some cases with relatively large number of
iterations (such as Ex.~6).  There is a case in which the GPGCD method
does not converge (Ex.~6).  Factors leading to such phenomena is under
investigation.

\begin{table*}[t]
  \centering
  \caption{Test results for $(m,d,n)$: $n$ input polynomials of degree
    $m$ with the degree of approximate GCD $d$. See
    Section~\ref{sec:exp} for details.} 
  \begin{tabular}{|c|c|c|c|c|c|c|c|c|c|}
    \hline
    Ex. & $(m,d,n)$ & \multicolumn{2}{c|}{\#Fail} & 
    \multicolumn{2}{c|}{Error}  & \multicolumn{2}{c|}{\#Iterations} &
    \multicolumn{2}{c|}{Time (sec.)} \\
    \cline{3-10}
    & & STLN & GPGCD & STLN & GPGCD & STLN & GPGCD & STLN & GPGCD \\
    \hline
    1 & $(10,5,3)$ & $0$ & $0$ & $2.31e\!-\!3$ & $2.38e\!-\!3$ &
    $5.50$ & $11.2$ & $1.17$ & $0.45$\\ 
    \hline
    2 & $(10,5,5)$ & $0$ & $0$ & $5.27e\!-\!3$ & $5.22e\!-\!3$&
    $4.70$ & $13.5$ & $3.10$ & $1.53$ \\ 
    \hline
    3 & $(10,5,10)$ & $0$ & $0$ & $5.48e\!-\!3$ & $5.62e\!-\!3$&
    $4.40$ & $17.9$ & $12.49$ & $8.59$ \\ 
    \hline
    4 & $(20,10,3)$ & $0$ & $0$ & $5.17e\!-\!3$ & $5.40e\!-\!3$ &
    $4.50$ & $12.0$ & $3.35$ & $1.52$ \\  
    \hline
    5 & $(20,10,5)$ & $0$ & $0$ & $5.89e\!-\!3$ & $5.85e\!-\!3$ &
    $4.40$ & $12.7$ & $10.37$ & $4.97$ \\ 
    \hline
    6 & $(20,10,10)$ & $0$ & $1$ & $6.31e\!-\!3$ & $6.20e\!-\!3$ &
    $4.00$ & $25.6$ & $44.62$ & $43.16$ \\ 
    \hline
    7 & $(40,20,3)$ & $0$ & $0$ & $5.32e\!-\!3$ & $5.39e\!-\!3$ &
    $4.90$ & $12.8$ & $13.60$ & $5.83$ \\ 
    \hline
    8 & $(40,20,5)$ & $0$ & $0$ & $6.01e\!-\!3$ & $5.97e\!-\!3$ &
    $4.30$ & $12.1$ & $41.46$ & $17.92$ \\ 
    \hline
    9 & $(40,20,10)$ & $0$ & $0$ & $6.41e\!-\!3$ & $6.25e\!-\!3$ &
    $4.10$ & $8.90$ & $200.88$ & $60.21$ \\ 
    \hline
  \end{tabular}
  \label{tab:appgcd}
\end{table*}

\section{Concluding Remarks}

Based on our previous research (\cite{ter2009}, \cite{ter2009b}), we
have extended our GPGCD method for more than two input polynomials
with the real coefficients.  We have shown that, at least
theoretically, our algorithm properly calculates an approximate GCD
under certain conditions for multiple polynomial inputs.

Our experiments have shown that, in the case that the number of
iterations is relatively small, the GPGCD method calculates an
approximate GCD efficiently with almost the same amount of
perturbations as the STLN-based method.  However, computing time of
the GPGCD method increases as the number of iterations becomes
larger; it suggests that we need to reduce the computing time of each
iteration in the GPGCD method for the cases with relatively large
number of iterations.  It is desirable to have more detailed
experiments for analyzing stability, performance for input polynomials
of larger degree, etc.

For the future research, generalizing this result to polynomials with
the complex coefficients will be among our next problems.  It is also
an interesting problem how the choice of $P_1$ affects the performance
of the algorithm.  Furthermore, one can also use arbitrary linear
combination to transform $\gcd(P_1,P_2,\ldots,P_n)$ to
$\gcd(P_1,a_2P_2+\cdots+a_nP_n)$.  This will reduce the size of the
generalized Sylvester matrix and will be another approach for
calculating approximate GCD.

\section*{Acknowledgments}

We thank Professor Erich Kaltofen for making their implementation for
computing approximate GCD available on the Internet.  We also thank
the anonymous reviewers for their valuable suggestions.

\bibliographystyle{splncs}

\def\cprime{$'$}

\end{document}